\documentclass{article}
\usepackage{graphicx}
\usepackage{lipsum}  
\usepackage{amsmath} 
\usepackage{enumitem} 
\usepackage{amsfonts} 
\usepackage{geometry} 
\usepackage{graphicx} 
\usepackage{caption} 
\usepackage{subcaption} 
\usepackage[ruled,vlined]{algorithm2e} 
\usepackage{listings} 
\usepackage{booktabs}

\begin{document}

\title{Fuzzy Hunter Optimizer: An Bio-Metaheuristic Algorithm Inspired
by L\'evy walks}
\author{
  {\small Hernández Rodríguez, Matías Ezequiel}\thanks{E-mail: \texttt{mhernandez@inidep.com}} \\ \\
  {\small Department of Mathematical Biology} \\
  {\small Instituto Nacional de Investigación y Desarrollo Pesquero} \\
  {\small Mar del Plata, Argentina} 
}

\date{} 

\maketitle

\begin{abstract}
This article introduces the Fuzzy Hunter Optimizer (FHO), a novel metaheuristic inspired by L\'evy diffuse visibility walk observed in predatory species and even in human behavior during the search for sustenance. To address a constrained optimization problem, we initialize a population of hunters in the search space. The hunter with the best fitness represents the food source. The other hunters move through the search space following a L\'evy walk. When they spot the food source, they move towards it, gradually abandoning the Levy walk. To model the hunters' visibility, we employ linear membership functions. In each iteration, the hunter with the best fitness becomes the food source. Unlike other metaheuristics, FHO parameters (visibility functions) do not require pre-calibration, since they adapt with iterations.
\end{abstract}


\section{Introduction}

Optimization represents a fundamental concept with widespread applications across diverse domains, ranging from engineering design and business strategy to internet routing and environmental sustainability. In the corporate sphere, enterprises aim to maximize profits while minimizing operational costs. Engineers are dedicated to optimizing the performance of their designs while simultaneously minimizing expenses. In the realm of sustainability, studies focus on minimizing the environmental impact associated with resource utilization. Virtually every facet of our activities is interconnected with the pursuit of optimization.
Consider, for instance, the process of vacation planning, where individuals seek to maximize their enjoyment while efficiently managing associated costs \cite{Diw, Gil, Gol, Hul, Len, Loz, Lupt, Set, Yan1}. Mathematical optimization serves as the primary tool to address these challenges. However, it is essential to acknowledge that certain algorithms, particularly gradient-based search methods, are inherently local in their search approach. These algorithms typically commence with initial assumptions and endeavor to enhance solution quality. In the case of unimodal functions, convexity guarantees that the final optimal solution also represents a global optimum. Conversely, in scenarios involving multimodal objectives, the search process may become ensnared in a local optimum.

An additional constraint arises when dealing with optimization problems featuring high-dimensional search spaces. Classical optimization algorithms often fall short in providing viable solutions due to the exponential growth of the search space as problem size increases, rendering exact methods impractical. Furthermore, the complexities inherent in real-world problems frequently hinder the ability to verify the uniqueness, existence, and convergence conditions that mathematical approaches demand \cite{Kel, Yan1}.

In contrast, metaheuristics have emerged as potent methodologies for tackling a wide array of intricate optimization challenges. They are particularly renowned for their capacity to uncover global optima, a trait that distinguishes them in the field. Prominent instances of metaheuristics include genetic algorithms (GA), ant colony optimization (ACO), and particle swarm optimization (PSO), among others. The versatility of metaheuristics extends their reach into the realms of science, technology, and engineering, as substantiated by numerous applications \cite{Hau, Mir1, Mir2, Yan1}.

Notably, recent years have witnessed a surge in interest surrounding nature-inspired, human behavior-mimicking, and physics-inspired metaheuristics. Examples of such innovative approaches encompass the Moth Search Algorithm (MSA), Grey Wolf Optimizer (GWO), Gold Rush Optimizer (GRO), Bat-Inspired Algorithm (BA), Ebola Optimization Search Algorithm (EOSA), the Gravitational Search Algorithm (GSA) and Tumoral and Tumoral Angiogenic Optimizer (TAO) \cite{Her2, Mir0, Oye, Ras, Yang, Zol}. These algorithms are tailored to address a spectrum of distinct optimization problems. However, it is imperative to underscore that no universally optimal algorithm exists, capable of delivering superior solutions across all conceivable optimization challenges. Certain algorithms exhibit enhanced performance in specific problem domains, accentuating the need for ongoing exploration and innovation in the realm of heuristic optimization algorithms.

This paper introduces a novel metaheuristic method called Fuzzy Hunter Optimizer (FHO), which employs concepts inspired by L\'evy's diffuse visibility walk, a navigational strategy observed in predatory species such as silky sharks, yellowfin tuna, swordfish, and even human beings during their quest for sustenance \cite{Bro, Gar, Rai, Sim}. In Section 2, we elaborate on the Levy flight. In Section 3, we expound on the metaheuristic proposed in this article, which is inspired by the Levy flight with fuzzy vision. In the same section, we introduce the concept of visibility. In Section 4, we assess the performance of our algorithm using benchmark functions. Additionally, we tackle classical engineering problems (cantilever beam design, pressure vessel design, and tension/compression spring). In all cases, we compare the performance with other algorithms. Finally, in Section 5, we present the conclusion.

\section{L\'evy walks}

L\'evy walks, also known as a L\'evy flight, is a remarkable stochastic process that has garnered significant attention in various scientific disciplines, including physics, ecology, and optimization. It characterizes a unique movement pattern observed in diverse animal species, including sharks, tuna, swordfish, and even humans, during their foraging activities.

L\'evy walks is characterized by the distribution of step lengths, where step lengths are not drawn from a conventional exponential or Gaussian distribution. Instead, they follow a power-law distribution, specifically a L\'evy stable distribution, with heavy tails. The probability density function (PDF) of a L\'evy stable distribution is given by:

\begin{equation} \label{equ1}
	p(s) \sim \frac{1}{|s|^{\beta}}, \mbox{ for } |s| > 0, 
\end{equation}
where, $s$ represents the step length, and $0<\beta \leq 3$ is the stability parameter, which determines the behavior of the distribution. The heavy tails of the distribution mean that occasionally, there are very long steps, which are interspersed with shorter steps.

The L\'evy walk exhibits two key characteristics:

\begin{enumerate}
\item	Long Jumps: Occasionally, an individual will take a particularly long step. These long steps can significantly influence the overall movement pattern. 
\item	Scale Invariance: The L\'evy walk is scale-invariant, meaning that it displays similar statistical properties at different spatial and temporal scales. This feature allows for efficient exploration of space, making it an adaptive strategy for foraging.
\end{enumerate}

L\'evy distribution, denoted as $p(s)$ in (\ref{equ1}), can be expressed, according to Mantegna's method, as follows:

\begin{equation} \label{equ2}
		p(s) = \frac{(\beta - 1)\Gamma(\beta - 1)\sin\left(\frac{\pi(\beta-1)}{2}\right)}{\pi s^{\beta}}, 		
\end{equation}
where $s$ is bigger than 0 and $\Gamma(x)$ is the gamma function.

In the next section we explain the metaheuristic proposed in this paper.

\section{Fuzzy Hunter Optimizer algorithm}

When some predators in the animal kingdom, such as sharks, and human hunter-gatherers, such as the Hadza of northern Tanzania, face challenges locating food sources, their movement patterns deviate from Brownian motion (which is analogous to the random movement observed in gas molecules) and adopt a strategy known as L\'evy walks, which was discussed above. Subsequently, upon detecting a potential food source, they exhibit a directional response while still incorporating an element of randomness into their movement.

Therefore, if we denote the positions of the food source and the hunter at time $t$ as $x_f(t)$ and $x(t),$ respectively, we may consider defining the hunter's position at the subsequent time step, $t+1,$ as follows:

\begin{equation}\label{equ3}
	x(t+1) =
	\begin{cases}
		x(t) + w r(x_f(t) - x(t)) & \text{if the hunter perceives the food source,} \\
		\alpha L(s) & \text{otherwise,}
	\end{cases}
\end{equation}
where, $L(s)$ is the step extracted from L\'evy walks using Mantegna's method, $w$ represents a constant, $\alpha$ denotes an acceleration factor, $r$ is a random vector uniformly distributed in the interval $[0, 1],$ and $L(s)$ is defined by equation (\ref{equ2}). Notably, the previous analysis did not incorporate the concept of visual detection of a food source. We proceed to formalize this concept through the introduction of a linear membership function, defined as follows:

\begin{equation}\label{equ4}
	v(x(t),x_f(t)) = \min\left(1, \max\left(\frac{\|x(t)-x_f(t)\|_2 - R_u}{R_v - R_u} , 0\right)\right).
\end{equation}

Here, $v(x(t),x_f(t))$ signifies the visibility of the food source, $x_f(t),$ as perceived by the hunter $x(t),$ subject to the conditions $R_v > R_u \geq 0.$ According to equation (\ref{equ4}), the hunter maintains full visibility ($v(x(t),x_f(t))=1$) when situated at a distance less than $R_u$ from the food source. Conversely, the visibility of the food source diminishes linearly to zero as the distance between the hunter and the food source reaches or exceeds $R_v.$ A hunter with zero visibility resorts to moving in accordance with a L\'evy walk strategy.

We define the hunter's position at time \(t+1\) as follows:

\begin{equation}\label{equ5}
    x(t+1)= x(t) + wrv(x(t),x_f(t))(x(t)-x_f(t)) + \left(1-v(x(t),x_f(t)\right) L(s).
\end{equation}

Equation (\ref{equ5}) models the behavior of the hunter concerning the food source. As the visibility $v(x(t),x_f(t))$ increases, the term $\left(1-v(x(t),x_f(t))\right)$ decreases, indicating the hunter's departure from the L\'evy walk to move toward the food source. One notable advantage is that visibility is dynamically adjusted in each iteration, eliminating the need for manual parameter calibration as required in some other metaheuristics.

In the context of minimizing a function $f(x)$ within a search space $\Omega\subseteq \mathbb{R}^n,$ the FHO metaheuristic initiates with a randomized population of hunters, representing potential solutions within $\Omega.$ The hunter exhibiting the best fitness is designated as the food source. At each iteration $t,$ hunter $i,$ $x_i(t),$ advances towards the food source, $x_f(t),$ following the equation:

\begin{equation}\label{equ6}
    x_i(t+1)= x_i(t) + wrv(x_i(t),x_f(t))(x_i(t)-x_f(t)) + \left(1-v(x_i(t),x_f(t)\right) L(s).
\end{equation}

In this article, we have chosen specific parameter values: $w=2,$ $\beta=0.8,$ $R_u=10^{-n}d(\Omega)$ and $R_v = 0.8R_u,$ where $d(\Omega)$ denotes the diameter of $\Omega$ and is defined as:

$$
d(\Omega) = \sup\left\{\|x - y\|_2:x,y\in \Omega\right\}.
$$
The calculation of the diameter of a set can often be accomplished with ease; however, in certain cases, a numerical algorithm may be necessary \cite{Har, Malan, Ram}. The selection of an appropriate numerical method is contingent upon the characteristics of $\Omega.$

The fundamental procedures of the FHO algorithm are succinctly encapsulated in the pseudocode illustrated in Algorithm \ref{alg_FHO}.

\begin{algorithm}
\caption{Fuzzy Hunter Optimizer}
\label{alg_FHO}  
\KwData{Objective function, $f(x);$ search space, $\Omega;$ dimension of problem, $n;$ parameter $w;$ population size, $N_{pop}$ and max number of iterations, $N_{max}$} 

\KwResult{Approximate optimum, $x_f,$ and $f(x_f)$}

Initialize the hunter population randomly inside the search region\;
Calculate the fitness of each search hunter\;
Calculate the diameter of $\Omega$\;
$x_f(0)\Leftarrow $ the best search hunter\;

\While{$t<N_{max}$}{
    \ForEach{$x(t)$ $\neq$ $x_f(t)$}{
        Update the position of the $x(t)$ from equation (\ref{equ6})\;
        \If{$f(x(t))<f(x_f(t))$}{
            		$x_f(t) \Leftarrow x(t)$\; 
            		}     
    }
\Return{$x_f$ and $f(x_f).$}
}

\end{algorithm}

\section{Validation and comparation}
In order to assess the performance of our algorithm, in subsection \ref{benfun}, we applied it to 10 benchmark functions and compared the results with others algorithms.
In subsection \ref{claengpro}, real-world problems are addressed (cantilever beam design, pressure vessel design, and tension/compression spring), and the results of our algorithm are compared with several well-known and new optimization algorithms.   

\subsection{Benchmark functions} \label{benfun}
Table \ref{Tab1} enumerates the reference functions employed in our experimental investigation.
To assess the effectiveness of FHO, we conducted a comparative analysis against both established and emerging algorithms, including the Harris Hawks optimizer (HHO), PSO, ant lion optimizer (ALO), whale optimization algorithm (WOA), Levy flight distribution (LFD), hunter-prey optimizer (HPO), and tunicate swarming algorithm (TSA). The performance metrics of these algorithms were originally documented in \cite{Naru}. In order to maintain congruity with the findings reported therein, the maximum iteration and population size are 500 and 30, respectively, and FHO is run 30 times independently.

\begin{table}[ht]
	\noindent\begin{tabular*}{\textwidth}{l @{\extracolsep{\fill}}cccc}
	\toprule
	Function& Dimension& Range& Minimum\\
	\hline
 	\footnotesize$f_1=\displaystyle\sum_{i=1}^{n}x_i^2$& 30& $[-100,100]^n$& 0\\
 	\footnotesize$f_2=\displaystyle\sum_{i=1}^{n} |x_i| + \displaystyle\prod_{i=1}^{n} |x_i|$& 30& $[-10,10]^n$& 0\\
 	\footnotesize$f_3=\displaystyle\sum_{i=1}^{n} \left(\displaystyle\sum_{j=1}^{i}x_j\right)^2$& 30& $[-100,100]^n$& 0\\
 	\footnotesize$f_4=\displaystyle\sum_{i=1}^{n}100\left(x_i-x_{i-1}^2\right)+(x_{i-1}-1)^2$& 30& $[-30,30]^n$& 0\\
 	\footnotesize$f_5=\max\{|x_i|, 1\leq i \leq n\}$& 30& $[-100,100]^n$& 0\\
 	\footnotesize$f_6=\displaystyle\sum_{i=1}^{n}-x_i\sin \sqrt{|x_i|}$& 30& $[-500,500]^n$& - 418.9829$\times n$\\
 	\footnotesize$f_7=10n+\displaystyle\sum_{i=1}^{n}\left[x_i^2-10\cos(2\pi x_i)\right]$& 30& $[-5.12,5.12]^n$& 0\\
 	\footnotesize$f_8=\frac{1}{4000}\displaystyle\sum_{i=1}^{n}x_i^2 - \displaystyle\prod_{i=1}^{n}\cos\left(\frac{x_i}{\sqrt{i}}\right)+1$& 30& $[-600,600]^n$& 0\\
 	\footnotesize$f_9=\displaystyle\sum_{i=1}^{n}(x_i+0.5)^2$& 30& $[-100,100]^n$& 0\\
 	\footnotesize$f_{10}=-20\exp\left(-0.2\sqrt{\frac{1}{n}\displaystyle\sum_{i=1}^{n}x_i^2}\right)-\exp\left(\frac{1}{n}\displaystyle\sum_{i=1}^{n}\cos(2 \pi x_i)\right)+20+e$& 30& $[-32,32]^n$& 0\\
	\hline
	\end{tabular*}
\caption{Classical benchmark functions.}
\label{Tab1}   
\end{table}

\begin{figure}[ht]
    \centering
    \includegraphics[width=0.8\textwidth]{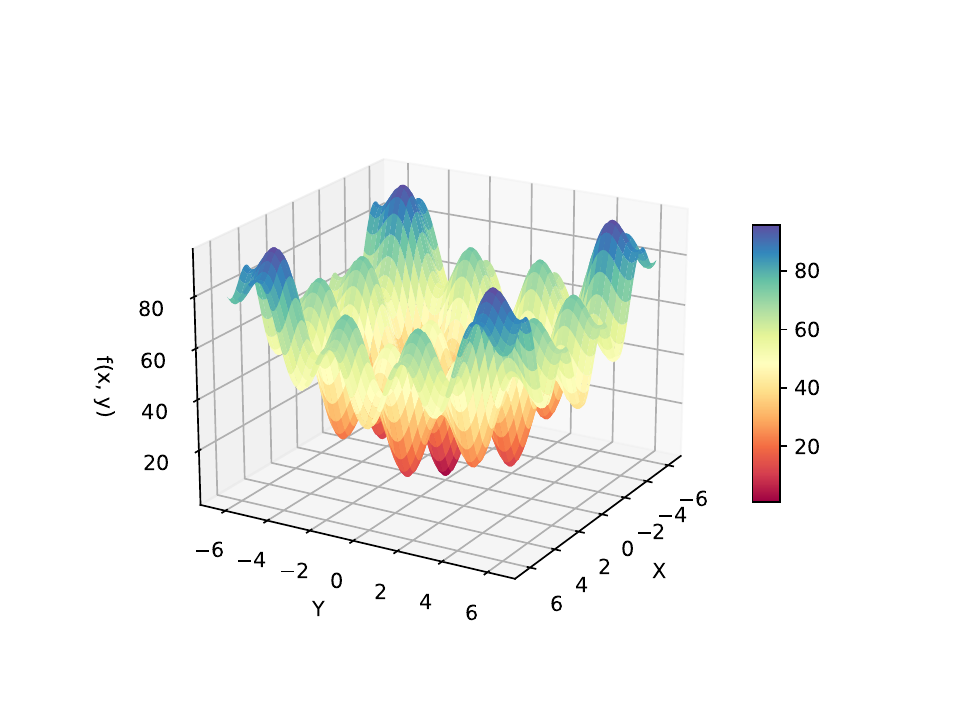}
    \caption{Eggcrate function: $x^2+y^2+25(\sin^2x+\sin^2y)$.}
    \label{figura2}
\end{figure}

\begin{figure}[ht]
\centering
  \begin{subfigure}[b]{0.4\textwidth}
    \includegraphics[width=\textwidth, height=\textwidth]{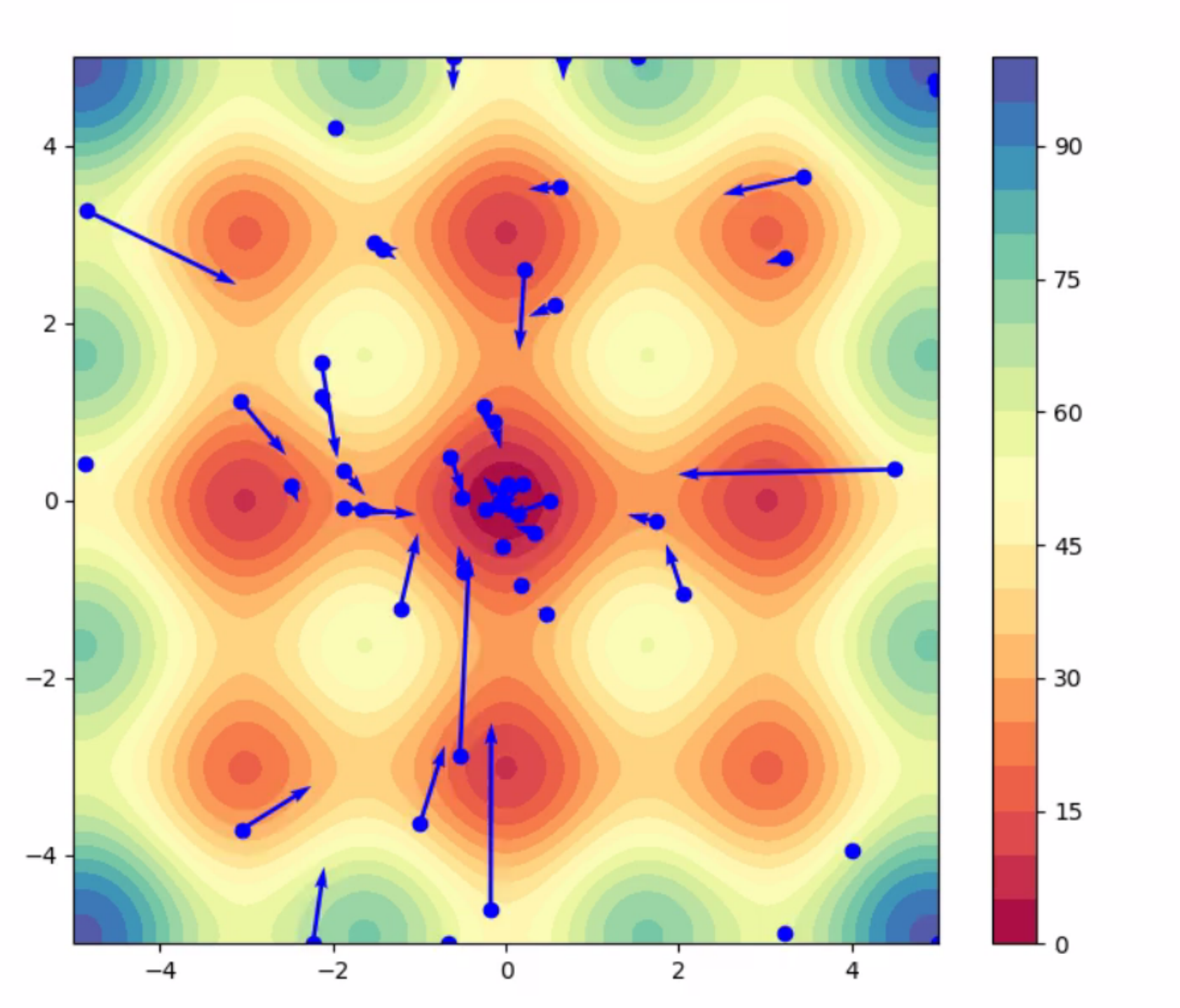}
    \caption{}
    \label{fig:f1}
  \end{subfigure}
  \begin{subfigure}[b]{0.4\textwidth}
    \includegraphics[width=\textwidth, height=\textwidth]{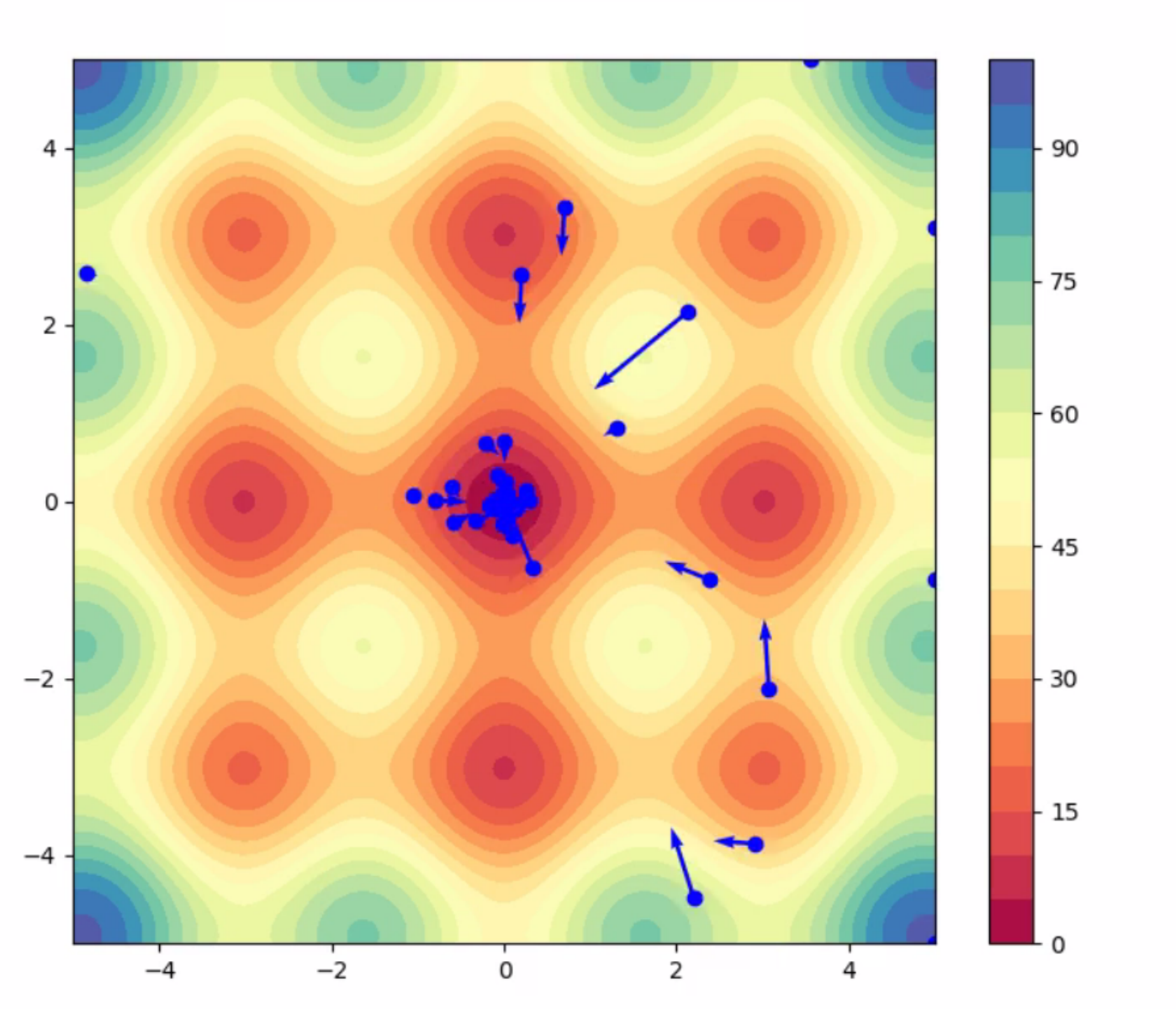}
    \caption{}
    \label{fig:f2}
  \end{subfigure}
  \begin{subfigure}[b]{0.4\textwidth}
    \includegraphics[width=\textwidth, height=\textwidth]{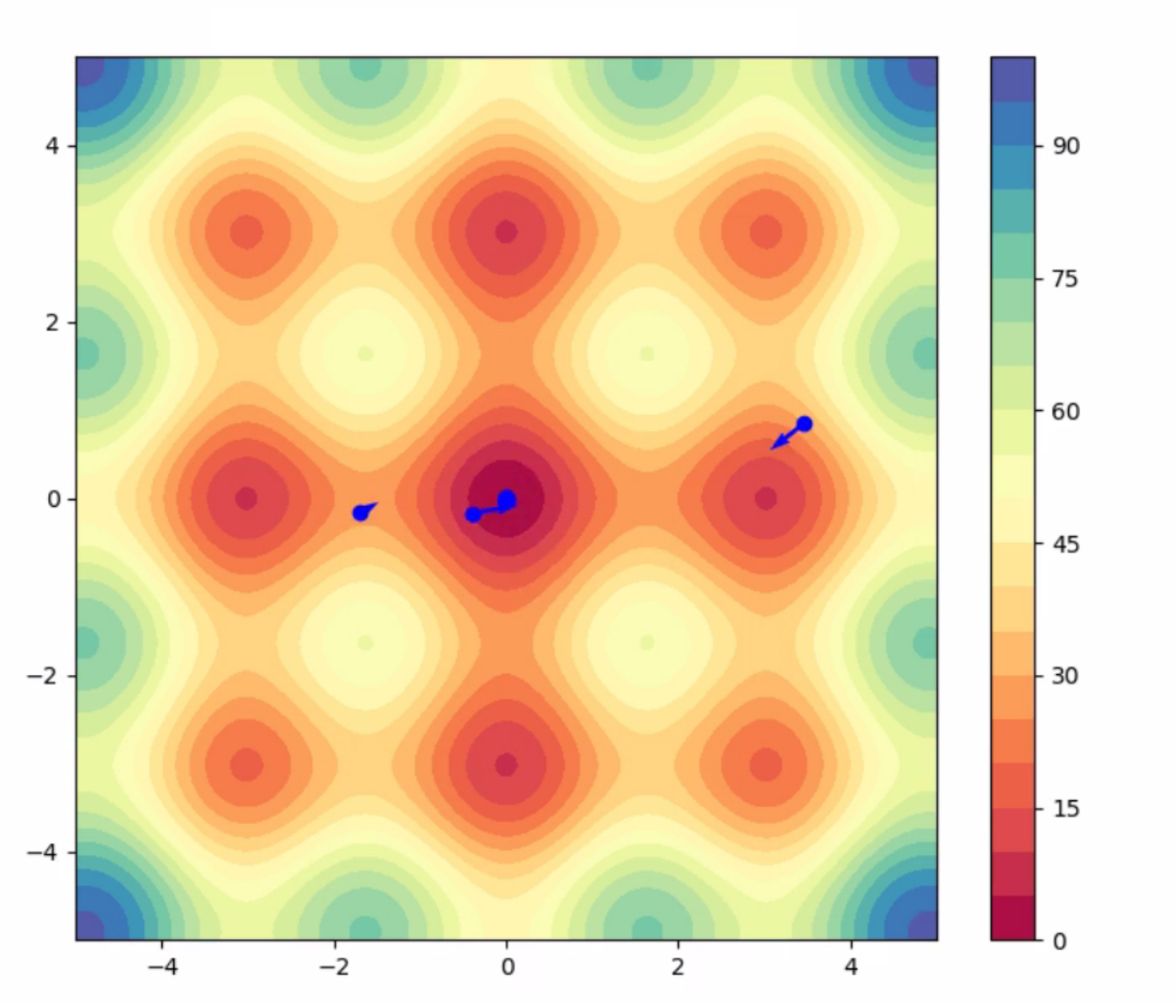}
    \caption{}
    \label{fig:f3}
  \end{subfigure}
  \begin{subfigure}[b]{0.4\textwidth}
    \includegraphics[width=\textwidth, height=\textwidth]{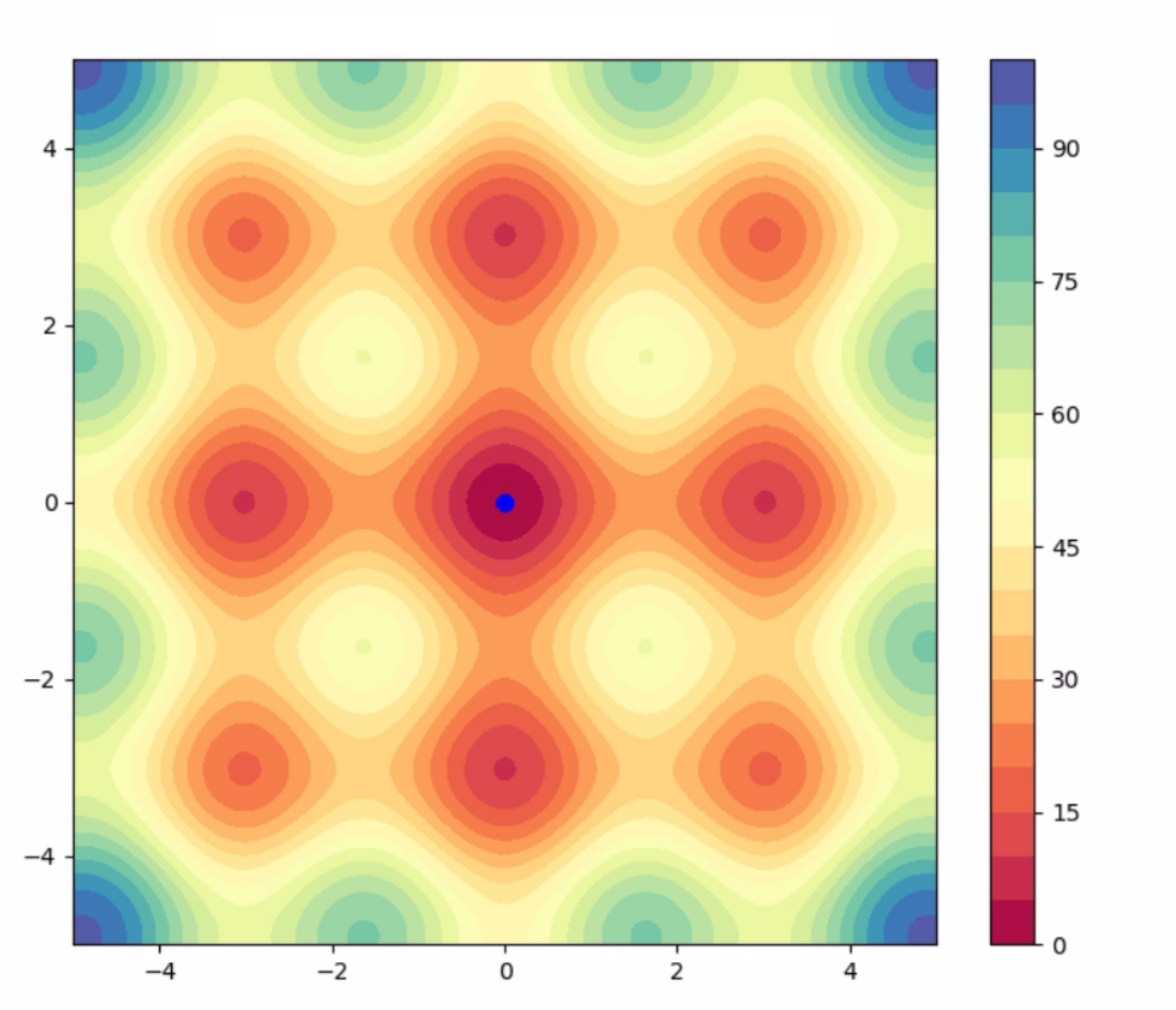}
    \caption{}
    \label{fig:f4}
  \end{subfigure}
  \caption{ Eggcrate function: Four instances depicting convergence in the FHO algorithm towards the minimum (a-d).}
  \label{figura1}
\end{figure}

\begin{table}[ht]
{\scriptsize
	\noindent\begin{tabular*}{\textwidth}{c @{\extracolsep{\fill}}ccccccccc}
	\toprule
	\small Function& & \scriptsize PSO& \scriptsize WOA& \scriptsize ALO& \scriptsize LFD& \scriptsize TSA& \scriptsize HHO& \scriptsize HPO& \scriptsize \textbf{FHO}\\
	\hline
		 & Min& 4.1929e-08& 9.1868e-86& 2.0758e-04& 1.5743e-07& 1.5625e-23& 9.6839e-117& 1.3483e-186& \textbf{7.6433e-37}\\
 	$f_1$& Mean& 2.7368e-06& 7.2644e-74& 0.0010& 3.1904e-07& 3.5086e-21& 1.1030e-93& 3.8960e-169& \textbf{1.1935e-31}\\
 		 & Std.& 7.8979e-06& 2.6139e-73& 8.3478e-04& 7.0385e-08& 7.5860e-21& 5.5816e-93& 0& \textbf{3.5983e-31}\\
 		 & Min& 8.7537e-05& 7.9595e-56& 0.8653& 2.2253e-04& 1.0225e-14& 4.0014e-60& 4.2018e-99& \textbf{8.2066e-17}\\
 	$f_2$& Mean& 0.0023& 1.2605e-50& 43.2560& 3.3647e-04& 1.1593e-13& 2.5192e-49& 2.9108e-92& \textbf{1.3564e-15}\\
 		 & Std.& 0.0042& 3.9289e-50& 47.8404& 6.0363e-05& 1.8668e-13& 1.3643e-48& 1.1177e-91& \textbf{1.3982e-15}\\
 		
 		 & Min& 48.6588& 1.5867e04& 657.0275& 5.0909e-07& 1.4390e-07& 6.0866e-102& 2.4279e-161& \textbf{2.8463e-08}\\
 	$f_3$& Mean& 236.4627& 4.6712e04& 4.4373e03& 1.3700e-06& 6.1660e-04& 1.8158e-71& 6.8323e-144& \textbf{9.7133e-05}\\
 		 & Std.& 147.1275& 1.5757e04& 2.5847e03& 5.0401e-07& 0.0013& 9.9440e-71& 3.7299e-143& \textbf{1.4946e-04}\\
 	
 		 & Min& 8.7916& 27.2605& 26.2370& 27.8536& 26.1806& 7.1629e-07& 22.8180& \textbf{7.4414e-05}\\
 	$f_4$& Mean& 48.5377& 28.0400& 330.3476& 28.0628& 28.3860& 0.0171& 23.7090& \textbf{9.4880}\\
 		 & Std.& 41.4386& 0.4365& 535.2123& 0.1291& 0.8620& 0.0223& 0.7436& \textbf{8.9630}\\

 		 & Min& 1.0616& 0.8494& 8.2454& 2.8538e-04& 0.0049& 6.0051e-57& 2.4541e-83& \textbf{3.8862e-02}\\
 	$f_5$& Mean& 2.6549& 47.0941& 16.6688& 3.4946e-04& 0.3257& 1.3491e-47& 5.4057e-77& \textbf{2.5377}\\
 		 & Std.& 0.9912& 28.9341& 5.2473& 4.5172e-05& 0.3510& 6.3104e-47& 1.6203e-76& \textbf{2.2139}\\	
	
 		 & Min& -7.8120e03& -1.2565e04& -8.5214e03& -5.0651e03& -6.9100e03& - 1.2569e04& -9.7233e03& \textbf{-1.2569e04}\\
 	$f_6$& Mean& -6.2641e03& -1.0020e04& -5.6298e03& -4.0512e03& -5.9093e03& -1.2569e04& -8.8439e03& \textbf{-1.2333e04}\\
 		 & Std.& 773.1355& 1.7649e03& 598.9912& 413.8118& 562.0738& 0.9315& 598.0833& \textbf{8.8631e02}\\		

 		 & Min& 21.8891& 0& 39.7995& 9.5044e-08& 115.9126& 0& 0& \textbf{2.2737e-13}\\
 	$f_7$& Mean& 46.5640& 1.8948e-15& 81.3609& 5.1798e-06& 184.6362& 0& 0& \textbf{1.3929e01}\\
 		 & Std.& 710.1715& 1.0378e-14& 24.3112& 5.4679e-06& 40.2045& 0& 0& \textbf{1.4891e01}\\	
 	
 		 & Min& 2.7242e-07& 0& 0.0167& 4.2184e-07& 0& 0& 0& \textbf{5.5511e-16}\\
 	$f_8$& Mean& 0.0231& 0& 0.0625& 8.5391e-07& 0.0053& 0& 0& \textbf{7.8456e-15}\\
 		 & Std.& 0.0281& 0& 0.0265& 2.2683e-07& 0.0074& 0& 0& \textbf{7.9339e-15}\\	
 	
 		 & Min& 6.8884e-08& 0.0846& 1.9036e-04& 0.8371& 2.8156& 7.7275e-07& 6.3266e-10& \textbf{4.4682e-31}\\
 	$f_9$& Mean& 1.5340e-06& 0.3742& 0.0012& 1.7965& 3.8706& 1.2635e-04& 1.2538e-07& \textbf{4.3235e-29}\\
 		 & Std.& 2.7959e-06& 0.2228& 7.4550e-04& 0.3084& 0.6509& 2.3701e-04& 4.1350e-07& \textbf{5.6880e-29}\\		
 		 & Min& 1.5447e-04& 8.8818e-16& 1.3414& 9.7034e-05& 2.0579e-12& 8.8818e-16& 8.8818e-16& \textbf{8.5709e-14}\\
 	$f_{10}$& Mean& 1.1543& 3.8488e-15& 4.3800& 1.2910e-04& 1.2947& 8.8818e-16& 8.8818e-16& \textbf{1.8054e00}\\
 		 & Std.& 0.8514& 2.8119e-15& 2.7934& 2.0618e-05& 1.6500& 0& 0& \textbf{1.0798e00}\\
	\hline
	\end{tabular*}
	}
\caption{Statistical results for minimizing the reference benchmark functions listed in Table \ref{Tab1}.}
\label{Tab2}   
\end{table}

The results presented in Table \ref{Tab2} demonstrate that the FHO algorithm performed well and is competitive with respect to the other analyzed algorithms. In most cases, it ranked among the top four performers. In some instances, such as the benchmark function $f_9$ it achieved the best performance.

\subsection{Classical engineering problems} \label{claengpro}

In this section, classic engineering problems are addressed (cantilever beam design, pressure vessel design, tension/compression spring and sustainable explotation renewable resource). The performance of our metaheuristic is compared with other algorithms, as reported in \cite{Naru}. The tackled problems are constrained optimization problems, which pose substantial challenges for metaheuristic algorithms.

Mathematically, constrained minimization problem are typically expressed as:

\begin{align} \label{eqOpt}
    & \min_{x\in\mathbb{R}^n} \quad  f(x)   \\
    & \text{Subject to:}  \nonumber \\
    & g_i(x)\leq 0, i =1,\dots, N \nonumber \\
    & h_j(x)=0, j=1,\dots, P \nonumber \\ 
    & l_h \leq x_h \leq u_h, h =1,\dots, n, \nonumber 
\end{align}
where $g_j$ are inequality constraints, $h_j$ are equality constraints, $l_i$ and $u_i$ are lower and upper bounds of $x_i,$ and $f(x)$ is the objective function that needs to be optimized subject to the constraints.

There are five prominent methodologies for handling constrained problems using metahuristics, which encompass:
1) \emph{Penalty functions}, 2) \emph{Specialized representations and operators}, 3) \emph{Repair algorithms}, 4) \emph{The separation of objectives and constraints}, and 5) \emph{Hybrid methods}.
Among these, penalty functions represent the most straightforward approach \cite{Coello}.
In this article, we utilize the FHO metaheuristic with the static penalty approach. Essentially, depending on the nature of the problem, as demonstrated in the application examples, we pursue two strategies.

One strategy, the more intuitive approach, involves transforming problem (\ref{eqOpt}) into the following problem:

\begin{align} \label{eqOpt1}
    & \min_{x\in\mathbb{R}^n} \quad  F(x)   \\
    & l_h \leq x_h \leq u_h, h =1,\dots, n \nonumber 
\end{align}
where $F(x)=\displaystyle\sum_{i=1}^Nr_i\max\{g_i(x),0\}+\displaystyle\sum_{j=1}^P c_j|h_j(x)|,$  functions $\max\{g_i(x),0\}$ and $|h_j(x)|$ measure the extent to which the constraints are violated, and $r_i$ and $c_j$ are known as penalty parameters.

In the other strategy, 

\begin{equation}\label{eqOpt2}
F(x) =
\begin{cases}
f(x) & \text{if the solution is feasible;} \\
K - \displaystyle\sum_{i=s} \frac{K}{m} & \text{otherwise,}
\end{cases}
\end{equation}
where $s$ is the number of constraints satisfied, $m$ is the total number of (equality and inequality) constraints and $K$ is a large constant ($K=1\times 10^9$) \cite{Mora}. 

\subsubsection{Cantilever Beam Design Problem} 

The Cantilever Beam Design Problem presents a noteworthy challenge within the domains of mechanics and civil engineering, primarily centered on the minimization of a cantilever beam's weight. In this problem, the beam comprises five hollow elements, each featuring a square cross-section. The objective is to ascertain the optimal dimensions for these elements while adhering to specific constraints (see Figure \ref{figura2}).

The mathematical expression governing this problem and its associated constraints are described by the following equation:

\begin{align} \label{equ7}
    & \min_{0.01\leq x_1,x_2,x_3,x_4,x_5\leq 100} \quad  f(x_1,x_2,x_3,x_4,x_5) = 0.06224(x_1+x_2+x_3+x_4+x_5)  \\
    & \text{Subject to:}  \nonumber \\ 
    & g(x):\frac{61}{x_1^3} + \frac{37}{x_2^3} + \frac{19}{x_3^3} + \frac{7}{x_4^3} + \frac{1}{x_5^3} \leq 1. \nonumber 
\end{align}

In order to address problem (\ref{equ7}), we reformulate it as an unconstrained problem as follows:

\begin{align} 
    & \min_{0.01\leq x_1,x_2,x_3,x_4,x_5\leq 100} \quad  F(x_1,x_2,x_3,x_4,x_5), 
\end{align}
where $F(x_1,x_2,x_3,x_4,x_5)=f(x_1,x_2,x_3,x_4,x_5)+\max\{g(x)-1, 0\}.$ 

Table \ref{Tab3} lists the best solutions obtained by FHO and various methods: artificial ecosystem-based optimization (AEO), ALO, coot optimization algorithm (COOT), cuckoo search algorithm (CS), gray prediction evolution algorithm based on accelerated even (GPEAae), HPO, interactive autodidactic school (IAS), multi-verse optimizer (MVO) and symbiotic organisms search (SOS) \cite{Naru}.

\begin{figure}[ht]
    \centering
    \includegraphics[width=0.8\textwidth]{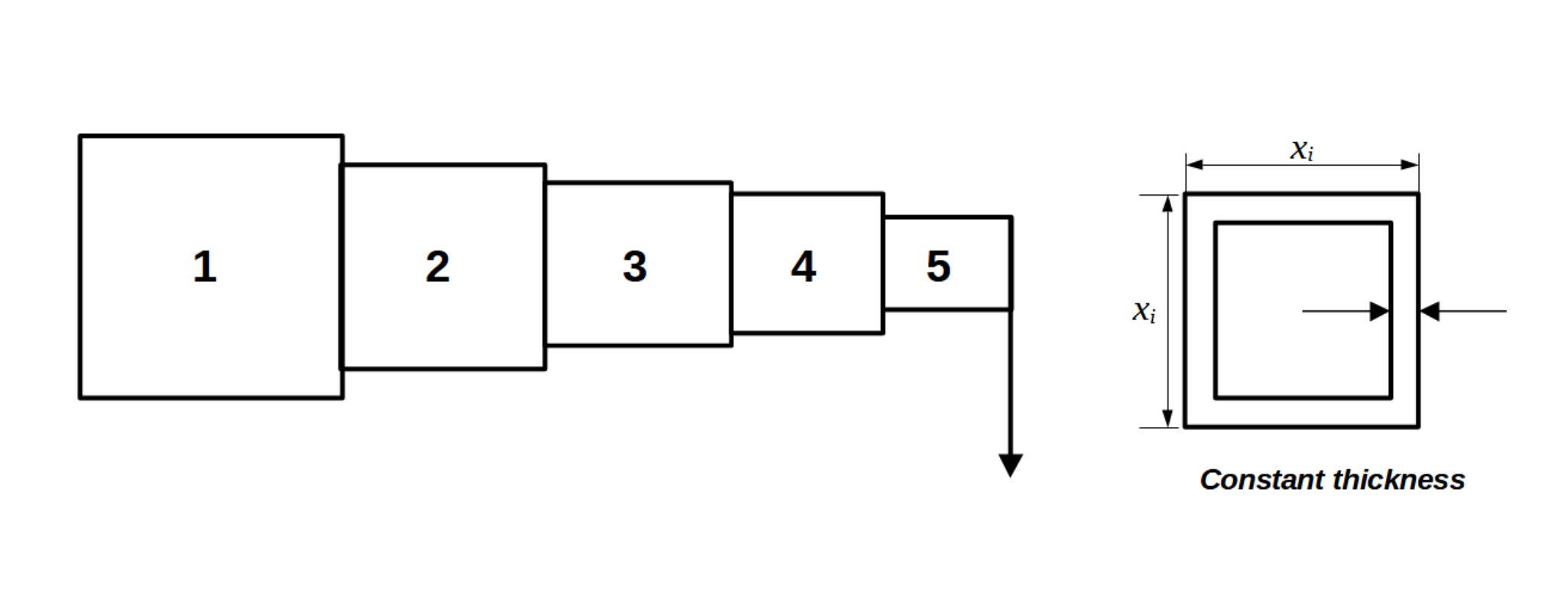}
    \caption{Cantilever Beam Design Problem.}
    \label{figura2}
\end{figure}

\begin{table}[ht]
{\small
    \noindent\begin{tabular*}{\textwidth}
    {l @{\extracolsep{\fill}}llllll }
    \toprule
    \multicolumn{1}{l}{\small Algorithm}& \multicolumn{5}{l}{\small Optimum variables}& \multicolumn{1}{c}{\small Optimum weight}\\
    \cline{2-6} 
        &$x_1$& $x_2$& $x_3$& $x_4$& $x_5$& \\
     \textbf{FHO}& \textbf{6.0421055}& \textbf{5.3377723}& \textbf{4.4720019}& \textbf{3.4819607}& \textbf{2.1409217}& \textbf{1.3365892} \\
	 AEO& 6.02885000& 5.31652100& 4.46264900& 3.50845500& 2.15776100& 1.33996500 \\   
     ALO& 6.01812000& 5.31142000& 4.48836000& 3.49751000& 2.15832900& 1.33995000 \\
    COOT& 6.02743657& 5.33857480& 4.49048670& 3.48343700& 2.13459100& 1.33657450 \\
       CS& 6.0089000& 5.30490000& 4.50230000& 3.50770000& 2.15040000& 1.33999000 \\
   GPEAae& 6.0148080& 5.30672400& 4.49323200& 3.50516800& 2.15378100& 1.33998200 \\
   HPO& 6.00552336& 5.30591367& 4.49474956& 3.51336235& 2.15423400& 1.33652825 \\ 
      IAS& 5.9914000& 5.30850000& 4.51190000& 3.50210000& 2.16010000& 1.34000000 \\
      MVO& 6.0239402& 5.30601120& 4.49501130& 3.49602200& 2.15272610& 1.33995950 \\
      SOS& 6.0187800& 5.30344000& 4.49587000& 3.49896000& 2.15564000& 1.33996000 \\
     \hline
    \end{tabular*}
    }
    \caption{Comparison results for the cantilever design problem.}
    \label{Tab3}   
\end{table}

\subsubsection{Pressure vessel design problem} 

\begin{figure}[ht]
    \centering
    \includegraphics[width=0.8\textwidth]{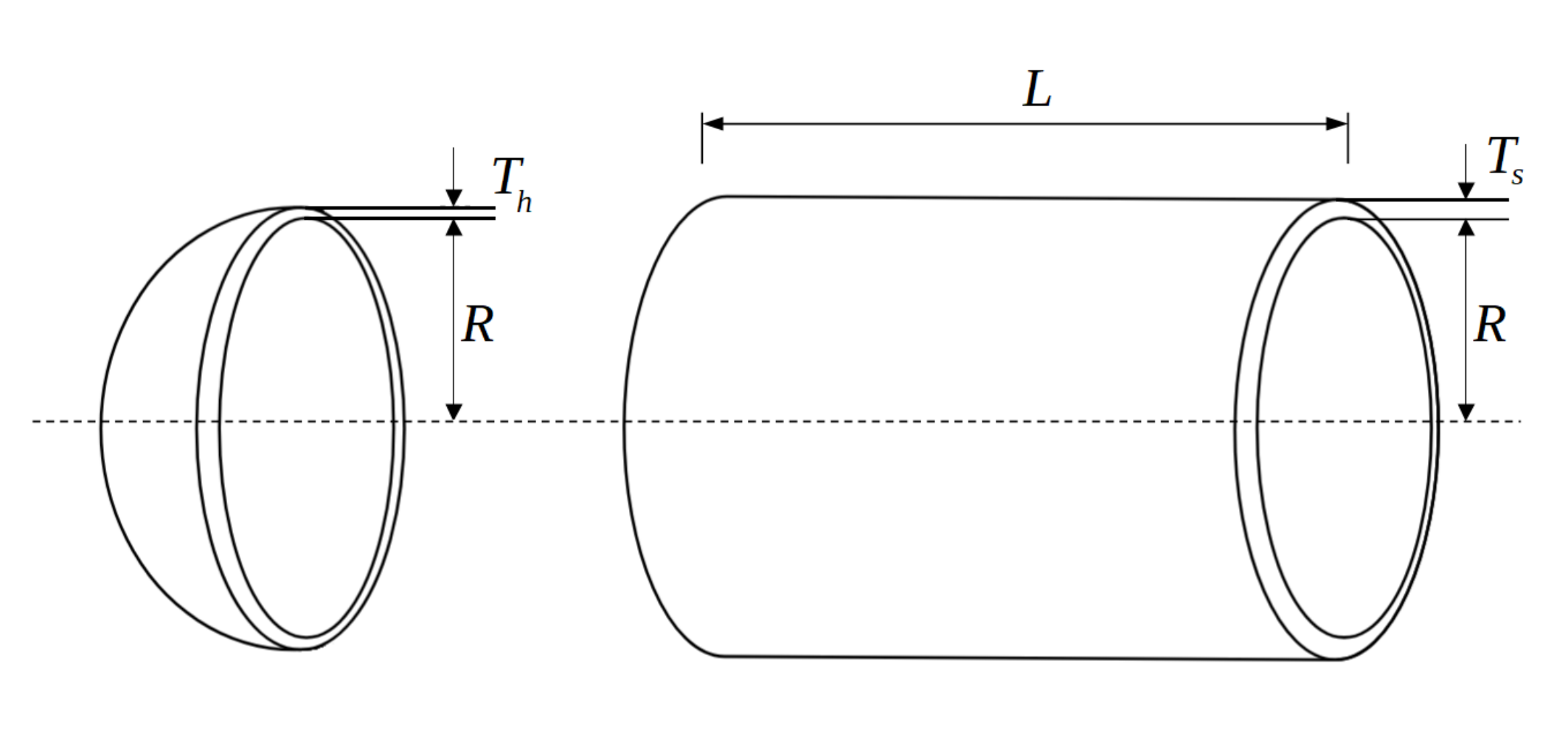}
    \caption{Pressure vessel design problem.}
    \label{figura3}
\end{figure}

The aim of this problem is to achieve cost minimization encompassing material, forming, and welding expenses associated with a cylindrical vessel, as illustrated in Figure \ref{figura3}. The vessel is enclosed at both ends with a hemispherical-shaped head. There are four key variables in this problem: $T_s$ ($x_1,$ thickness of the shell), $T_h$ ($x_2,$ thickness of the head), $R$ ($x_3,$ inner radius) and $L$ ($x_4,$ length of the cylindrical section of the vessel, not including the head).

This problem is governed by four constraints. The formulation of the problem is articulated as follows \cite{Kan}:

\begin{align} \label{equ8}
    & \min \quad  f(x) = 0.6244x_1x_3x_4 + 1.7781x_2x_3^2 +3.1661x_1^2x_4 + 19.84x_1^2x_3  \\
    & \text{Subject to:}  \nonumber \\ 
    & g_1(x): -x_1 + 0.0193x_3\leq 0 \nonumber \\
    & g_2(x): -x_2 + 0.00954x_3\leq 0 \nonumber \\
    & g_3(x): -\pi x_3^2x_4 - \frac{4}{3}\pi x_3^3 +1296000 \leq 0 \nonumber \\ 
    & g_4(x): x_4 - 240\leq 0 \nonumber \\
    & 0 \leq x_1 \leq 99 \nonumber \\
    & 0 \leq x_2 \leq 99 \nonumber \\
    & 10 \leq x_3 \leq 200 \nonumber \\
    & 10 \leq x_4 \leq 200. \nonumber 
\end{align}

To solve problem (\ref{equ8}), we transform it into an unconstrained problem following strategy (\ref{eqOpt2}), with the number of constraints being $m=4.$ 

Comparative results with other algorithms are presented in Table \ref{Tab4}. Algorithms under consideration include AEO, BA, charged system search (CSS), GA, GPEAae, Gaussian quantum-behaved particle swarm optimization (G-QPSO), GWO, HHO, HPO, hybrid particle swarm optimization (HPSO), mothflame optimization (MFO), sine cosine gray wolf optimizer (SC-GWO), water evaporation optimization(WEO) and WOA. 

\begin{table}[ht]
{\small
    \noindent\begin{tabular*}{\textwidth}{l @{\extracolsep{\fill}}lllll}
    \toprule
    \multicolumn{1}{l}{\small Algorithm}& \multicolumn{4}{l}{\small Optimum variables}& \multicolumn{1}{l}{\small Optimum cost}\\
    \cline{2-5} 
        &$T_s$& $T_h$& $R$& $L$& \\
      \textbf{FHO}& \textbf{0.8375030}& \textbf{0.4139782}& \textbf{43.3939372}& \textbf{161.2185336}& \textbf{5994.6845509} \\
      AEO& 0.8374205& 0.413937& 43.389597& 161.268592& 5994.50695 \\   
      BA& 0.812500& 0.437500& 42.098445& 176.636595& 6059.7143 \\
      CSS& 0.812500& 0.437500& 42.103624& 176.572656& 6059.0888 \\
      GA& 0.812500& 0.437500& 42.097398& 176.654050& 6059.9463 \\
      GPEAae& 0.812500& 0.437500& 42.098497& 176.635954& 6059.708025 \\ 
      G-QPSO& 0.812500& 0.437500& 42.0984& 176.6372& 6059.7208 \\
      GWO& 0.8125& 0.4345& 42.089181& 176.758731& 6051.5639 \\
      HHO& 0.81758383& 0.4072927& 42.09174576& 176.7196352& 6000.46259 \\
      HPO& 0.778168& 0.384649& 40.3196187& 200& 5885.33277 \\
      HPSO& 0.812500& 0.437500& 42.0984& 176.6366& 6059.7143 \\
      MFO& 0.8125& 0.4375& 42.098445& 176.636596& 6059.7143 \\
      SC-GWO& 0.8125& 0.4375& 42.0984& 176.63706& 6059.7179 \\
      WEO& 0.812500& 0.437500& 42.098444& 176.636622& 6059.71 \\
      WOA& 0.812500& 0.437500& 42.0982699& 176.638998& 6059.7410 \\
     \hline
    \end{tabular*}
    }
    \caption{Comparison of results for the pressure vessel design problem.}
    \label{Tab4}   
\end{table}

\subsubsection{Tension/compression spring design problem}
\begin{figure}[ht]
    \centering
    \includegraphics[width=0.6\textwidth]{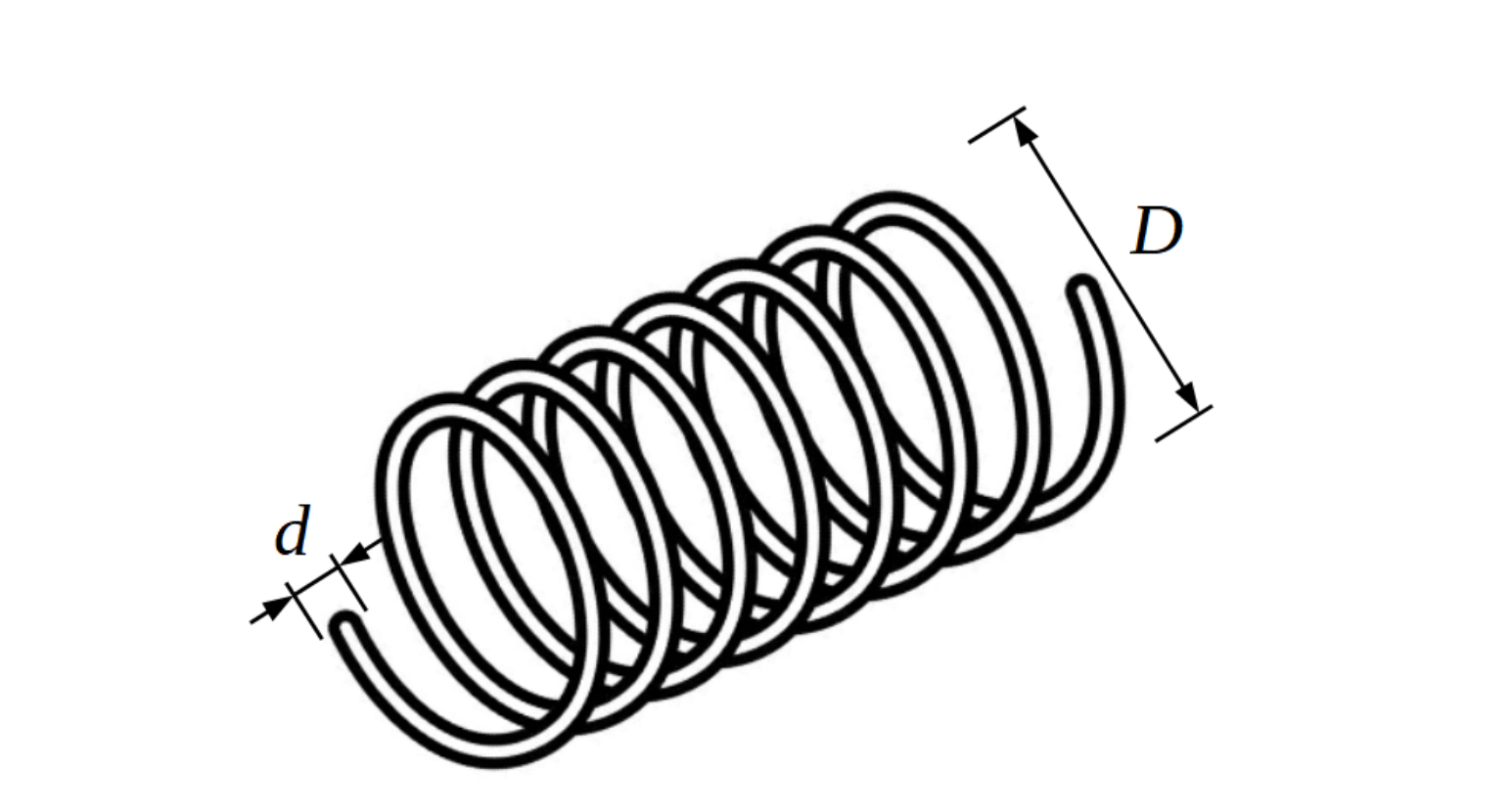}
    \caption{Tension/compression spring design.}
    \label{figura4}
\end{figure}

The engineering test problem utilized in this study revolves around optimizing the design of tension/compression springs. The aim of this problem is to minimize the weight of a tension/compression spring characterized by three crucial parameters: the number of active loops ($N$), the average coil diameter ($D$), and the wire diameter ($d$). The visual representation of the spring's geometric details and its associated parameters can be found in Figure \ref{figura4}.

The spring design problem is subject to a set of inequality constraints, formally defined in Equation (\ref{equ9}).

\begin{align} \label{equ9}
    & \min \quad  f(x) = (x_3+2)x_2x_1  \\
    & \text{Subject to:}  \nonumber \\ 
    & g_1(x): 1 - \frac{x_2^3x_3}{71785x_1^4}\leq 0 \nonumber \\
    & g_2(x): \frac{4x_2^2-x_1x_2}{12566(x_2x_1^3-x_1^4)} + \frac{1}{5108x_1^2} - 1 \leq 0 \nonumber \\
    & g_3(x): 1 - \frac{140.45x_1}{x_2^2x_3} \leq 0 \nonumber \\ 
    & g_4(x): \frac{x_2+x_1}{1.5} - 1 \leq 0 \nonumber \\
    & 0.05 \leq x_1 \leq 2 \nonumber \\
    & 0.25 \leq x_2 \leq 1.3 \nonumber \\
    & 2.00 \leq x_2 \leq 15, \nonumber 
\end{align}
where $x_1 = d,$ $x_2 = D$ and $x_3=N.$

To solve problem (\ref{equ9}), we have applied the methodology described in approach (\ref{eqOpt2}), with the number of constraints in this problem being $m=4.$ Table \ref{Tab5} presents a comparative analysis of the results obtained using the FHO algorithm alongside other heuristic optimization methods in the context of this specific problem. This table provides comprehensive information concerning the values of each design variable and their respective optimization outcomes. 
Algorithms considered in the comparison of results are the following: AEO, BA, COOT, CPSO, GPEAae, GWO, HPO, HHO, stochastic fractal search (SFS), spotted hyena optimizer (SHO), salp swarm algorithm (SSA), water cycle algorithm (WCA), WEO and WOA. 

\begin{table}[ht]
{\small
    \noindent\begin{tabular*}{\textwidth}{l @{\extracolsep{\fill}}llll}
    \toprule
    \multicolumn{1}{l}{\small Algorithm}& \multicolumn{3}{l}{\small Optimum variables}& \multicolumn{1}{l}{\small Optimum weigth}\\
    \cline{2-4} 
        &$N$& $D$& $d$&  \\
       \textbf{FHO}& \textbf{9.4875998}& \textbf{0.3919440}& \textbf{0.0531127}& \textbf{0.0127014} \\
       AEO& 10.879842& 0.361751& 0.051897& 0.0126662 \\   
       BA& 11.2885& 0.35673& 0.05169& 0.012665 \\
       COOT& 11.34038& 0.35584& 0.05165& 0.012665293 \\
       CPSO& 11.244543& 0.357644& 0.051728& 0.012674 \\
       GPEAae& 11.294000& 0.356631& 0.051685& 0.012665 \\
       GWO& 11.28885& 0.356737& 0.05169& 0.012666 \\
       HHO& 11.138859& 0.359305355& 0.051796393& 0.012665443 \\
       HPO& 11.21536452893& 0.3579796674& 0.0517414615& 0.012665282823723 \\
       SFS& 11.288966& 0.356717736& 0.051689061& 0.012665233 \\
	   SHO& 12.09550& 0.343751& 0.051144& 0.012674000 \\
	   SSA& 12.004032& 0.345215& 0.051207& 0.0126763 \\
       WCA& 11.30041& 0.356522& 0.05168& 0.012665 \\
       WEO& 11.294103& 0.356630& 0.051685& 0.012665 \\
       WOA& 12.004032& 0.345215& 0.051207& 0.0126763 \\
     \hline
    \end{tabular*}
    }
    \caption{Comparison of results for the Tension/compression spring problem.}
    \label{Tab5}   
\end{table}

\newpage
\section{Conclusions}

In this paper, we present a new metaheuristic algorithm called Fuzzy Hunter Optimizer (FHO). This algorithm is inspired by the L\'evy walks that many hunters follow in the animal kingdom. We have modeled the visibility of hunters using linear membership functions, allowing us to use learning parameters that adjust over time without having to calibrate them beforehand. According to the results obtained (see tables \ref{Tab1}-\ref{Tab5}), our algorithm has shown good performance in addressing the problem of minimizing 10 benchmark functions and classical engineering problems, which involve constrained optimization problems.

For future work, the proposed algorithm can be used in different fields of study. Furthermore, we aim to develop a version of our algorithm adapted to address multi-objective optimization problems.

\end{document}